\documentclass[12pt]{article}
\usepackage{amsmath,amssymb,amsfonts,amscd}
\usepackage[cp1251]{inputenc}

\newtheorem{thm}{Theorem}[section]
\newtheorem{cor}[thm]{Corollary}
\newtheorem{lem}[thm]{Lemma}
\newtheorem{prop}[thm]{Proposition}

\title {Resolutions of free partially commutative monoids}
\author {L.\,~Polyakova}

\begin {document}
\date{}
\maketitle

\begin{quote}
{\small {\bf Key words:} monoid homology, free partially
commutative monoid, free resolution, homological dimension.

{\bf Abstract.} A free resolution of free partially commutative
monoids is constructed and with its help the homological dimension
of these monoids is calculated. }
\end{quote}

A free partially commutative monoid is a monoid generated by
elements some of which commute (see the precise definition in
Section 1). Homology of these monoids appeared in an article
\cite{HusTk} by A.\,~Husainov and V.\,~Tkachenko in connection
with constructing the homology groups of asynchronous transition
systems. In \cite{HusTezisy} Husainov proposes the following
Conjecture:

\medskip

{\bf Conjecture.} {\it Let $\Sigma$ be a finite set and $M$ be a
free partially commutative monoid whose generating set is
$\Sigma$. If there are no distinct letters $a_1,a_2,\dots
a_{n+1}\in \Sigma$ such that $a_ia_j=a_ja_i$ for every $1\leq
i<j\leq n+1$, then the monoid $M$ has the homological dimension
$\leq n$.}

\medskip

In this paper we construct a free resolution for a free partially
commutative monoid and with its help prove the Husainov's
Conjecture. We follow the ideas of D.\,~Cohen who built in
\cite{Cohen} a resolution for the so-called graph product of
groups, given resolutions for factors. The presentation of the
graph product with the help of direct and free amalgamated
products played the leading role at that. However the additional
difficulties appear while using this method for monoids.

Section 1 is devoted to basic definitions and facts concerned with
free partially commutative monoids and free amalgamated products.
In Section 2 the desired resolution is constructed. If the
opposite is not specified all considered modules are right.

\section{Preliminaries}

In the subsequent text we follow mainly \cite{DiekMet} in
considering free partially commutative monoids and
\cite{CliffPrOrig1} in considering monoid free amalgamated
products.

Let $\Sigma$ be a finite set called the alphabet. We denote by
$\Sigma^{\ast}$, and its elements are called the words, the free
monoid generated by  $\Sigma$. The notation $alph(x)$ is used for
all letters of $\Sigma$ that appear in a word $x\in
\Sigma^{\ast}$.

Let $I\subseteq \Sigma\times \Sigma$ be a symmetric irreflexive
binary relation over the alphabet $\Sigma$ called the commutation
relation. The complement of $I$ is denoted by
$D=\Sigma\times\Sigma\setminus I$.

A monoid $M(\Sigma,I)$, which has  a presentation $<\Sigma|
\{ab=ba, (a,b)\in I\}>$, is called a free partially commutative
monoid.

An undirected graph without loops $\Gamma(M)$ can be uniquely
compared to the free partially commutative monoid $M(\Sigma, I)$
in the following way: the vertex set of $\Gamma(M)$ is $\Sigma$,
and the edges connect commuting vertices.

Important tools to work with free partially commutative monoids
are the  Projection Lemma and Levi's Lemma (\cite{DiekMet}). To
formulate the former we need a definition. Let $A\subseteq \Sigma$
and $I_A=(A\times A)\cap I$ be an induced commutation relation.
{\it The projection} is a homomorphism $\pi_A:M(\Sigma,
I)\rightarrow M(A,I_A)$ which erases all the letters from a word
which do not belong to $A$. In other words, for $a\in \Sigma$ we
have

\[ \pi_A(a)=\left\{
  \begin{array}{ll}
    a, & a\in A, \\
    1, & a\not\in A.
  \end{array}
  \right.
\]

If $I_A=\emptyset$, then $\pi_A$ is a projection of $M(\Sigma,I)$
onto the free monoid $A^\ast$. If $A=\{a,b\}$, we write
$\pi_{a,b}$ instead of  $\pi_{\{a,b\}}$.

\begin{lem}\label{LemProject}{\bf (Projection Lemma, \cite{DiekMet})}
Elements $u,v\in M(\Sigma,I)$ are equal if and only if
$\pi_{a,b}(u)=\pi_{a,b}(v)$ for all  $(a,b)\in D$.
\end{lem}

\begin{lem}\label{LemLevi}{\bf (Levi, \cite{DiekMet})}
Let $t,u,v,w\in M(\Sigma,I)$. The following assertions are
equivalent:
\begin{description}
  \item[1)] $tu=vw$;
  \item[2)] there exist $p,q,r,s\in M(\Sigma,I)$ such that $t=pr$,
  $u=sq$, $v=ps$, $w=rq$ with $rs=sr$ and $alph(r)\cap alph(s)=\emptyset$.
\end{description}
\end{lem}

To write down the elements of the free partially commutative
monoid the so called Foata normal form is used. It is defined in
the following way. Let $\Sigma$ be totally ordered. A word $x\in
M(\Sigma, I)$ is in the Foata normal form if either it is the
empty word or  there exist an integer $n>0$ and non-empty words
$x_i$ $(1\leq i\leq n)$, such that
\begin{description}
  \item[1)] $x=x_1x_2 \dots  x_n$;
  \item[2)] for each $i$ the word $x_i$ is a product of distinct pairwise commuting
  letters, the letters of $x_i$ being written with regard to
  ordering introduced on $\Sigma$;
  \item[3)] for each $1\leq i<n$ and for each letter $a$ of $x_{i+1}$ there exists a letter $b$ of $x_i$,
  such that   $(a,b)\in D$.
\end{description}

The following theorem  holds:

\begin{thm}\label{ThmNormForm}({\rm \cite{DiekMet}, \cite{Lallem79Orig}})
Every element of $M(\Sigma,I)$ has a unique Foata normal form.
\end{thm}

For example, the Foata normal form of the $n$-th power $a^n$
$(a\in \Sigma)$ consists of $n$ factors equal $a$.

We recall briefly the concepts of the monoid free product and the
free amalgamated product.

A monoid free product $\Pi^\ast\{M_j, j\in J\}$ or  simple
$\Pi^\ast M_j$ is built for a monoid family $\{M_j, j\in J\}$
provided $M_i\cap M_j={1}$, $i\ne j$ (see, for instance,
\cite[V.2., $\S 9.4$]{CliffPrOrig1}). It consists if the
single-element sequence (1) and  all nonempty sequences
$(a_1,\dots ,a_k)$, such that $a_j\ne 1$, $a_j\in M_{i(j)}$,
$j=1,\dots ,k$ and $i(j)\ne i(j+1)$, $j=1,\dots ,k-1$.

For each $j\in J$ a canonical isomorphic embedding
$\chi_{j}:M_j\rightarrow \Pi^\ast M_j$ can be defined as follows:
$\chi_j(a)=a$; and $M_j$ can be identified with its canonical
image. Hence, we can suppose that $\Pi^\ast M_j$ is generated by
its submonoids $M_j$. The element $(a_1, \dots ,a_k )\in \Pi^\ast
M_j$ can be written as $a_1\dots a_k$.

A monoid free amalgamated product is built for a family $[\{M_j,
j_J\}; U; \{\varphi_j, j_J \}]$ which is called a monoid amalgam.
Here $\{ M_j, j\in J\}$ and $U$ are monoids. We assume again that
$M_i\cap M_j = 1$, $i\ne j$, and in the free product $\Pi^\ast
M_j$ each monoid $M_j$ is identified with its canonical image.
Homomorphisms $\varphi_j, j\in J$ are the embeddings of the monoid
$U$ into the monoids $M_j$, such that the unit of $U$ is mapped
onto the common unit of the monoids $M_j$.

We define the relation $\nu$ over $\Pi^\ast M_j$ assuming
\[\nu =\{ (u_i,u_j)| u_i=\varphi_i(u), u_j=\varphi_j(u),\, \mbox{for some}\, i,j\in J, u\in U \}.\]
Let $\sim_\nu$ be the least congruence containing $\nu$. A monoid
$\Pi^\ast \{M_j, j\in J \}/\sim_\nu$ is called a free product of
the amalgam $[\{M_j, j_J\}; U; \{\varphi_j, j_J \}]$ or a free
amalgamated product and is denoted by $\Pi^\ast_U M_j$.

The free amalgamated product can be described in terms of
generators and defining relations, namely, the following
proposition holds:

\begin{prop}\label{PropSvAmProizv}({\rm  \cite{CliffPrOrig1}})
Let $[M_j; U; \varphi_j]$ be a monoid amalgam and the monoid $U$
has a presentation $<Y|\pi>$ where $\pi\subset Y^\ast \times
Y^\ast$. Then there exist such sets $X_j$ that $Y\subseteq X_j$,
$X_i\cap X_j=Y$, if $i\ne j$, and there exist such relations
$\sigma_j\subset X_j^\ast \times X_j^\ast$ that $M_j=<X_j,
\sigma_j>$, $\sim_\pi = \sim_{\sigma_j}\cap (Y^\ast \times
Y^\ast)$ for each $j$. If $X=\bigcup_{j\in J}X_j$, $\sigma =
\bigcup_{j\in J}\sigma_j$, then the free amalgamated product
$\Pi^\ast_U M_j$ has a presentation $<X|\sigma>$.
\end{prop}

\section{Constructing resolutions}

We assume that the tensor product is considered over the ring
$\mathbb{Z}$ if it is not specified. Also, for the monoid $M$ we
write  ``$M$-module'' instead of ``$\mathbb{Z}M$-module'' and
$\otimes_M$ means that the tensor product is considered over the
ring $\mathbb{Z}M$.

First, we discuss how a resolution for a free commutative monoid
looks like. The free commutative monoid $M$ with $n$ generators
$a_1,a_2, \dots, a_n$ is a direct product of $n$ infinite cyclic
monoids $M^1,M^2,\dots , M^n$ with generators $a_1,a_2, \dots a_n$
respectively. For each of them the resolution looks  like
\[ 0\rightarrow [a_j]\mathbb{Z}M^j \stackrel{\partial^j}\rightarrow \mathbb{Z}M^j \stackrel{\varepsilon^j}\rightarrow \mathbb{Z} \rightarrow 0\]
where $[a_j]\mathbb{Z}M^j$ is a free $M^j$-module with one
generator $[a^j]$ and $\partial^j[a^j]=a^j-1$.

Denote by $X^j$ the complex $0\rightarrow [a^j]\mathbb{Z}M^j
\rightarrow \mathbb{Z}M^j$. Reasoning similarly for monoids as in
\cite[IV, \S 6]{Macl63Orig} we obtain the resolution for $M$ as
the tensor product of complexes $X^j$ with the augmentation
$\varepsilon = \varepsilon^1\otimes\dots \otimes\varepsilon^n$:

\begin{equation} 0\rightarrow X_n \stackrel{\delta_n}\rightarrow \dots
\rightarrow X_2 \stackrel{\delta_2}\rightarrow
X_1\stackrel{\delta_1}\rightarrow X_0=\mathbb{Z}M
\stackrel{\varepsilon}\rightarrow  \mathbb{Z} \rightarrow 0,
\label{rez1}\tag{$\ast$}
\end{equation}
where
\[ X_k = \sum_{m_1+\dots + m_k=n}X^1_{m_1}\otimes \dots \otimes X^n_{m_n}.\]
The number of the summands in this sum equals  $C_n^k$,  and each
of them can be identified with the free $M$-module with one
generator $[a_{i_1}a_{i_2}\dots a_{i_k}]$ where $1\leq
i_1<i_2<\dots i_k\leq n$ and
$X^{i_j}_{m_{i_j}}=[a_j]\mathbb{Z}M^{i_j}$, $j=1,2,\dots, k$.

The boundary homomorphisms $\delta_k$ are of the form:
\[\delta_k[a_{i_1}\dots a_{i_k}]=\sum_{j=1}^k[a_{i_1}\dots \widehat{a_{i_j}}\dots
a_{i_k}](a_{i_j}-1)(-1)^{j-1}.
\]

Before we turn to the main Theorem we prove two lemmas.

\begin{lem}\label{LemSvPodmod}
Let $M=M(\Sigma,I)$ be a free partially commutative monoid,
$\Sigma_0\subset \Sigma$, $I_0=(\Sigma_0\times\Sigma_0)\cap I$ and
$M_0=M(\Sigma_0,I_0)$. Then the monoid ring $\mathbb{Z}M$ is a
free (left) $\mathbb{Z}M_0$-module.
\end{lem}

{\bf Proof.} To make sure that $M_0$ is really a submonoid of $M$
we build an embedding $i:M_0\rightarrow M$. Let $[x]_{M_0}$ be an
element of $M_0$ defined by $x\in \Sigma_0^\ast$, and $[y]_M$ be
an element of $M$ defined by $y\in \Sigma^\ast$. Now set
$i([x]_{M_0})= [x]_M$.

This mapping is correct and is a homomorphism, for instance, in
view of \cite[V.1, \S 1.12, corollary 1.29]{CliffPrOrig1}. We
prove that $i$ is an injection. If $x,y\in \Sigma^\ast_0$ and
$i([x]_{M_0})=[x]_M=[y]_M=i([y]_{M_0})$, then $x$ can be obtained
from $y$ by successive transpositions of neighboring letters
$a,b$, such that $(a,b)\in I$. However $(a,b)\in I_0$ since
$x,y\in \Sigma_0^\ast$. Hence, $[x]_{M_0}=[y]_{M_0}$.

We order the elements of $\Sigma$ in such a way that all the
element of $\Sigma_0$ precede all the elements of
$\Sigma\setminus\Sigma_0$.

To prove that $\mathbb{Z}M$ is a free module we construct its
basis. Consider the set $B$ which consists of the monoid unit and
all the elements $t\in M$, such that the presentation $t$ in the
Foata normal form $t=w_1w_2\dots w_n$ has the following property:
$w_1$ consists only of those letters which belong to
$\Sigma\setminus\Sigma_0$.

Notice that if $u\in B$, then it cannot be presented in the form
$u=su_0$ where $s\in M_0\setminus 1$. Indeed, suppose such a
presentation exists. Consider the presentation of $u$ in the Foata
normal form $u=u_1u_2 \dots u_n$ and the letter $x$ which is the
first letter of the word $s$. The first occurrence of $x$ in the
word $u_1u_2 \dots u_n$ belongs to some $u_j$, $j>1$. Then by the
definition of the Foata normal form there exists such a letter $y$
of $u_{j-1}$ that $(x,y)\in D$. Two projections $\pi_{x,y}(u_1u_2
\dots u_n)$ and $\pi_{x,y}(su_0)$ do not coincide in the free
monoid $\{x,y\}^\ast$, since the first letter of $\pi_{x,y}(u_1u_2
\dots u_n)$ is $y$, but the first letter of $\pi_{x,y}(su_0)$ is
$x$, that contradicts  Lemma~\ref{LemProject}.

Show that each element $w\in M$ can be presented in the form
$w=au$ where $a\in M_0$, $u\in B$. To find such a presentation it
is sufficient to consider the following procedure. Present $w$ in
the Foata normal form: $w=w_1w_2\dots w_n$. If the letters of the
alphabet $\Sigma_0$ do not occur in the word $w_1$, then $w=1\cdot
w$ where $w\in B$. Otherwise, $w_1=a_1u_1$, where, in view of the
introduced on $\Sigma$  order, the words $a_1$ and $u_1$ can be
chosen in such a way that $a_1\in M_0$ and $u_1$ does not contain
the letters of the alphabet $\Sigma_0$. Consider the word
$u_1w_2\dots w_n$, present it in the Foata normal form
$u_1w_2\dots w_n=w_1^1w_2^1\dots w_n^1$ and again ``separate'' the
element of $M_0$ in the word $w_1^1$:  $w_1^1= a_2u_2$. Continuing
similarly we obtain finally the decomposition of the form
$w=a_1\dots a_k u_kw_1^kw_2^k\dots w_{n_k}^k = au$ where $a=a_1
\dots  a_n\in M_0$ and $u= u_k w_1^k w_2^k\dots w_{n_k}^k$ is an
element of the set $B$.

Thus we know how to decompose the elements of the monoid $M$ by
the elements of $B$. Now obviously we can decompose the elements
of the ring $\mathbb{Z}M$ by the elements of $B$ with the
coefficients in the ring $\mathbb{Z}M_0$.

To show that such a decomposition is unique it is sufficient to
prove that the element  $w$ of the monoid $M$ does not have two
different decompositions. Suppose $w=au= bv$ where $a,b\in M_0$,
$\{u,v\}\in B$. Then Lemma~\ref{LemLevi} implies that there exist
$p,q,r,s\in M$, such that $a=pr$, $b=ps$, $u=sq$, $v=rq$. From
$a=pr$, $b=ps$ it follows particulary that $r,s\in M_0$. But then
the presentations $u=sq$ and $v=rq$ contradict the fact $u,v\in B$
if $r=s=1$ does not hold. In this  case $a=b=p$, $u=v=q$ and the
decompositions $a\cdot u= b\cdot v$ coincide.

Therefore, the set $B$ is a basis that proves the lemma.
\rule{1ex}{1ex}

{\bf Remark.} Note that Lemma~\ref{LemSvPodmod} always holds for
groups, i.e. if $G$ is an arbitrary group and $G_0$ is its
subgroup, then the group ring $\mathbb{Z}G$ is a free
$\mathbb{Z}G_0$-module (see., for instance, \cite[I, \S
3]{Brown82Orig}). At the same time, if we choose a submonoid of
$M(\Sigma, I)$ which is not  free partially commutative, the
statement generally does not hold.

{\bf Example.} Let $M=<a>$ be an infinite cyclic monoid and
\linebreak $M_0=\{1,a^2,a^3,\dots ,a^n,\dots\}$ be a submonoid of
$M$. Then $\mathbb{Z}M$ is not a free $M_0$-module. Suppose the
contrary, then there exists a basis $B$ in $\mathbb{Z}M$ and each
of the elements of $\mathbb{Z}M$ has the unique decomposition by
this basis over the ring $\mathbb{Z}M_0$. Consider the elements
$1,a$. None of them can be presented as $xu$ where $x\in
M_0\setminus 1$, $u\in M$. Hence, both of them are contained in
$B$. But then the element $a^3=a^2\cdot a=a^3\cdot 1$ has two
decompositions by the elements of the basis.

Notice that the ring $\mathbb{Z}$ becomes a trivial $M_0$-module
if we set $n\cdot a =n$ for each $a\in M_0$, $n\in \mathbb{Z}$.
Then the tensor product $\mathbb{Z}\otimes_{M_0}\mathbb{Z}M$
exists and is an $M$-module, since $\mathbb{Z}M$ is an $M$-module.
This remark allows us to formulate the following Lemma.

\begin{lem}\label{LemTochnStroka}
Let $M$ be  a monoid, $M_0,M_1,M_2$ be its submonoids and
$M=M_1\ast_{M_0}M_2$. Then the following sequence of $M$-modules
\[0\rightarrow \mathbb{Z}\otimes_{M_0}\mathbb{Z}M \stackrel{i}\rightarrow \mathbb{Z}\otimes_{M_1}\mathbb{Z}M \oplus
\mathbb{Z}\otimes_{M_2}\mathbb{Z}M \stackrel{p}\rightarrow
\mathbb{Z} \rightarrow 0.
\]
is exact.
\end{lem}
{\bf Proof.} Define $i$ and $p$ in the following way. For
arbitrary $w,u,v\in M$ we set $i(1\otimes_{M_0} w)=1\otimes_{M_1}
w + 1\otimes_{M_2} w$, $p(1\otimes_{M_1} v)=1$, $p(1\otimes_{M_2}
u)=-1$. It is easy to check that in terms $\mathbb{Z}$ and
$\mathbb{Z}\otimes_{M_0} \mathbb{Z}M$ the sequence is exact. We
prove the exactness in the middle term. Since $pi(1\otimes_{M_0}
\sum_{w\in\Sigma} l_w w)=\sum_{w\in\Sigma}l_w-
\sum_{w\in\Sigma}l_w=0$, then ${\rm Im}\,i\subseteq {\rm Ker}\,p$.

The inverse involving must be proved. An element
$1\otimes_{M_1}\sum n_u u + 1\otimes_{M_2}\sum m_v v$ belongs to
${\rm Ker}\,p$ if and only if $\sum n_u =\sum m_v$. These sums can
be presented in the form:
\[ \sum n_u=n_+ - n_-,\,\,\,  \sum m_v=m_+ - m_-
\]
where $n_{+}, m_+$ are the sums of all the positive coefficients
and $n_-, m_-$ the sums of all the negative ones respectively. If
$n_+=m_+$, then $n_-=m_-$. It means that the sum
$1\otimes_{M_1}\sum n_u u$ contains  the same quantity of the
``plus''-sign summands of the form $1\otimes_{M_1} a$ ($a\in M$)
as the sum $1\otimes_{M_2}\sum m_v v$. Similarly for the
``minus''-sign summands. Suppose  $n_+\ne m_+$, for instance $n_+>
m_+$. Then we add  $n_+ - m_+$ summands of the form
$1\otimes_{M_2}b$ for some $b\in M$ to the sum $1\otimes_{M_2}\sum
m_v v$ and subtract them. Then the numbers of ``plus''-sign
summands and ``minus''-sign summands coincide for both of the
sums, since $m_+ +(n_+ - m_+)-(m_- + (n_+ - m_+))=n_+ - n_-$.

The following step is to construct a preimage under  the action of
$i$ for the element of the form  $1\otimes_{M_1}u +
1\otimes_{M_2}v$ where $u,v\in M$. After solving this problem we
will be able to find a preimage for all the elements   of ${\rm
Ker}\, p$ due to the reasoning stated above, i.e. to show that
${\rm Ker}\,\subseteq {\rm Im}\,i$ and finish the proof of the
Lemma.

To find the inverse image we present $u$ and $v$ in the form:
\[ u=a_1b_1a_2b_2\dots a_lb_l;\,\, v=c_1d_1c_2d_2\dots c_sd_s
\]
where $a_i,c_j\in M_1$, $b_i, d_j\in M_2$, $i=1,2,\dots, l$,
$j=1,2,\dots, s$. Then the element
\[
\begin{array}{l}
1\otimes_{M_0}w =
1\otimes_{M_0}[\displaystyle\sum_{k=1}^{l-1}(b_ka_{k+1}b_{k+1}\dots
a_lb_l - a_{k+1}b_{k+1}\dots a_lb_l) +b_l + \\
+\displaystyle\sum_{j=1}^{s}(c_jd_{j}c_{j+1}d_{j+1}\dots c_sd_s -
d_{j}c_{j+1}d_{j+1}\dots c_sd_s )]
\end{array}
\]
is a preimage of $1\otimes_{M_1}u + 1\otimes_{M_2}v$ under the
action of $i$. Indeed
\[
\begin{array}{l}
i(1\otimes_{M_0}w) = 1\otimes_{M_1}[b_1a_2b_2\dots a_lb_l +
\displaystyle\sum_{k=1}^{l-1}(-a_{k+1}b_{k+1}\dots a_lb_l +
b_{k+1}a_{k+2}b_{k+2}\dots a_lb_l) + \\
+ \displaystyle\sum_{j=1}^{s}(c_jd_{j}c_{j+1}d_{j+1}\dots c_sd_s -
d_{j}c_{j+1}d_{j+1}\dots c_sd_s )]+ \\
+1\otimes_{M_2}[\displaystyle\sum_{k=1}^{l-1}(b_ka_{k+1}b_{k+1}\dots
a_lb_l -
a_{k+1}b_{k+1}\dots a_lb_l) +b_l + \\
+ c_1d_1c_2d_2\dots c_sd_s +
\displaystyle\sum_{j=1}^{s-1}(-d_{j}c_{j+1}d_{j+1}\dots
c_sd_s + c_{j+1}d_{j+1}\dots c_sd_s ) - d_s]=\\
=1\otimes_{M_1}b_1a_2b_2\dots a_lb_l + 1\otimes_{M_2}(b_l+c_1d_1c_2d_2\dots c_sd_s-d_s)=\\
=1\otimes_{M_1}a_1b_1a_2b_2\dots a_lb_l +
1\otimes_{M_2}c_1d_1c_2d_2\dots c_sd_s
=1\otimes_{M_1}u+1\otimes_{M_2}v,
\end{array}
\]
what required. \rule{1ex}{1ex}

Now we are ready to prove the main theorem. Let $M(\Sigma, I)$ be
a free partially commutative monoid with a totally ordered
generating set $\Sigma$. Let  $\Gamma(M)$ be its graph and $r_k$
be the number of complete subgraphs with $k$ vertices in graph
$\Gamma(M)$. Let $F_k$ be a free $M$-module with  $r_k$
generators. We  denote each of such generators  as  $[a_1\dots
a_k]$ by putting the ascending sequence of vertices of
corresponding subgraph. We denote the $M$-module homomorphisms
$\delta_k:F_k\rightarrow F_{k-1}$, $k>1$ by setting for generators
\[ \delta_k[a_1\dots a_k]= \sum_{j=1}^k[a_1\dots \widehat{a_j}\dots
a_k](a_j-1)(-1)^{j-1}.
\]
Besides we set  $\delta_1[a]=a-1$ and $\varepsilon(a)=1$ for all
$a\in \Sigma$, and thus define homomorphisms
$\delta_1:F_1\rightarrow \mathbb{Z}M$ and
$\varepsilon:\mathbb{Z}M\rightarrow \mathbb{Z}$. The following
theorem holds:
\begin{thm}\label{ThmRezolv}
The sequence of  $M$-modules  and their homomorphisms

\begin{equation}
\dots\rightarrow F_n \stackrel{\delta_n}\rightarrow \dots
\rightarrow F_2 \stackrel{\delta_2}\rightarrow
F_1\stackrel{\delta_1}\rightarrow F_0=\mathbb{Z}M
\stackrel{\varepsilon}\rightarrow  \mathbb{Z} \rightarrow 0
\label{rez2}\tag{$\ast\ast$}
\end{equation}
is a free resolution of module $\mathbb{Z}$ over $\mathbb{Z}M$.
\end{thm}
{\bf  Proof.} It has to be  proved only the exactness of this
sequence. We use the induction on the number of generators of $M$.

If $M$ is a free commutative monoid (i.e. its graph is complete),
then the resolution (\ref{rez1}) coincides with a sequence
(\ref{rez2}) for this monoid, since the generators $[a_{i_1}\dots
a_{i_k}]$, $1\leq i_1<\dots < i_k\leq n$ of modules $X_k$ are in
one-one  correspondence with complete subgraphs  with the vertices
$a_{i_1},\dots ,a_{i_k}$. Particulary, in case when $M$ is
generated only by one generator, we obtain the induction
assumption.

Suppose $\Gamma(M)$ is not complete.  Then there exist two
vertices $x$ and $y$ which are not adjacent. Consider the
subgraphs $\Gamma_0=\Gamma\setminus\{x,y\}$,
$\Gamma_1=\Gamma\setminus x$, $\Gamma_2=\Gamma\setminus y$ and
their corresponding submonoids of the monoid $M$: $M_0(\Sigma_0,
I_0)$, $M_1(\Sigma_1, I_1)$ and $M_2(\Sigma_2, I_2)$. We have for
them:
\[
\begin{array}{l}
\Sigma_0=\Sigma\setminus\{x,y\};\,\,\Sigma_1=\Sigma\setminus
x;\,\,\Sigma_2=\Sigma\setminus y;\\
I_j=(\Sigma_j\times\Sigma_j)\cap I,\,j=0,1,2.
\end{array}
\]
Since for monoids $M_0, M_1, M_2$ the relations $\Sigma_1\cap
\Sigma_2=\Sigma_0$; $I_0=I_1\cap (\Sigma_0\times \Sigma_0)=
I_2\cap (\Sigma_0\times \Sigma_0)$ hold, then by
Proposition~\ref{PropSvAmProizv}  the free amalgamated product
$M_1\ast_{M_0}M_2$ has a presentation $<\Sigma_1\cup  \Sigma_2| \{
ab=ba, (a,b)\in I_1\cup I_2\}> =<\Sigma, \{ ab=ba, (a,b)\in I\}>$,
i.e. it coincides with the monoid $M$.

Further we apply the induction assumption to $M_0, M_1, M_2$. Let
\[ \dots \rightarrow
F_2^j \stackrel{\delta_2^j}\rightarrow
F_1^j\stackrel{\delta_1^j}\rightarrow \mathbb{Z}M_j
\stackrel{\varepsilon^j}\rightarrow  \mathbb{Z} \rightarrow 0,\,
j=0,1,2
\]
be the resolutions for these monoids. We consider their tensor
product with $\mathbb{Z}M$ over $\mathbb{Z}M_j$, $j=0,1,2$
respectively. By Lemma~\ref{LemSvPodmod} $\mathbb{Z}M$ is a free
$M_j$-module $(j=0,1,2)$, thus the functors  $\otimes_{M_j}
\mathbb{Z}M$ are exact, therefore, the sequences remain exact, the
modules $F_k^j\otimes_{M_j}\mathbb{Z}M$ being free $M$-modules.
Further consider the commutative diagram consisting of free
$M$-modules:

\begin{picture}(330,240)(-150,-50)
\put(-89,170){$\vdots$}\put(50,170){$\vdots$}\put(159,170){$\vdots$}
\put(-88,160){\vector(0,-1){15}}\put(51,160){\vector(0,-1){15}}\put(160,160){\vector(0,-1){15}}

\put(-140,130){0}\put(-132,135){\vector(1,0){15}}\put(-110,130){$F_2^0\otimes_{M_0}\mathbb{Z}M$}
\put(-45,135){\vector(1,0){15}}
\put(-20,130){$F_2^1\otimes_{M_1}\mathbb{Z}M \oplus
F_2^2\otimes_{M_2}\mathbb{Z}M$} \put(130,135){\vector(1,0){15}}
\put(155,130){$F_2$}\put(170,135){\vector(1,0){15}}\put(190,130){0}
\put(-42,138){$\scriptstyle{i_2}$}
\put(133,138){$\scriptstyle{p_2}$}

\put(-88,122){\vector(0,-1){15}}\put(51,122){\vector(0,-1){15}}\put(160,122){\vector(0,-1){15}}

\put(-140,90){0}\put(-132,95){\vector(1,0){15}}\put(-110,90){$F_1^0\otimes_{M_0}\mathbb{Z}M$}
\put(-45,95){\vector(1,0){15}}
\put(-20,90){$F_1^1\otimes_{M_1}\mathbb{Z}M \oplus
F_1^2\otimes_{M_2}\mathbb{Z}M$} \put(130,95){\vector(1,0){15}}
\put(155,90){$F_1$}\put(170,95){\vector(1,0){15}}\put(190,90){0}
\put(-42,98){$\scriptstyle{i_1}$}
\put(133,98){$\scriptstyle{p_1}$}

\put(-88,82){\vector(0,-1){15}}\put(51,82){\vector(0,-1){15}}\put(160,82){\vector(0,-1){15}}

\put(-145,50){0}\put(-137,55){\vector(1,0){15}}\put(-121,50){$\mathbb{Z}M_0\otimes_{M_0}\mathbb{Z}M$}
\put(-47,55){\vector(1,0){15}}
\put(-31,50){$\mathbb{Z}M_1\otimes_{M_1}\mathbb{Z}M \oplus
\mathbb{Z}M_2\otimes_{M_2}\mathbb{Z}M$}
\put(133,55){\vector(1,0){15}}
\put(150,50){$\mathbb{Z}M$}\put(173,55){\vector(1,0){15}}\put(193,50){0}
\put(-43,58){$\scriptstyle{i_0}$}
\put(135,58){$\scriptstyle{p_0}$}

\put(-88,42){\vector(0,-1){15}}\put(51,42){\vector(0,-1){15}}\put(160,42){\vector(0,-1){15}}

\put(-140,10){0}\put(-132,15){\vector(1,0){15}}\put(-105,10){$\mathbb{Z}\otimes_{M_0}\mathbb{Z}M$}
\put(-42,15){\vector(1,0){15}}
\put(-15,10){$\mathbb{Z}\otimes_{M_1}\mathbb{Z}M \oplus
\mathbb{Z}\otimes_{M_2}\mathbb{Z}M$}
\put(130,15){\vector(1,0){15}}
\put(155,10){$\mathbb{Z}$}\put(170,15){\vector(1,0){15}}\put(190,10){0}
\put(-38,18){$\scriptstyle{i}$} \put(134,18){$\scriptstyle{p}$}

\put(-88,2){\vector(0,-1){15}}\put(51,2){\vector(0,-1){15}}\put(160,2){\vector(0,-1){15}}

\put(-91,-30){0}\put(48,-30){0}\put(157,-30){0}

\end{picture}

As it has been already noticed the left and the middle columns are
exact. The bottom row is exact by Lemma~\ref{LemTochnStroka}. The
second (from the bottom) row is a sequence of modules
\[
0\rightarrow \mathbb{Z}M \stackrel{i_0}\rightarrow
\mathbb{Z}M\oplus \mathbb{Z}M \stackrel{p_0}\rightarrow
\mathbb{Z}M \rightarrow 0
\]
where  $i_0(a)=a\oplus a$ and $p_0(b\oplus c)=b-c$ for all
$a,b,c\in M$. From here it can be easy shown that it is also
exact.

To prove the exactness of other rows we note that each complete
subgraph, which is contained in $\Gamma_0$, is contained in
$\Gamma_1$ and in $\Gamma_2$ simultaneously, and each complete
subgraph, which is contained in $\Gamma$, is contained either in
$\Gamma_1$ or in $\Gamma_2$, since the vertices $x$ and $y$ are
not adjacent.

Let elements $[c_1],\dots ,[c_l]$ be generators of $F_n^0$;
$[a_1],\dots ,[a_m],[c'_1],\dots ,[c'_l]$ be generators of
$F_n^1$; $[b_1],\dots ,[b_k],[c''_1],\dots ,[c''_l]$ be generators
of  $F_n^2$ (for the sake of simplicity we denote the generators
with a single letter and give them accents for not to mistake what
direct summand they belong to). Then the generators of $F_n$ are
$[a_1],\dots ,[a_m],[b_1],\dots ,[b_k],[c_1],\dots ,[c_l]$. The
homomorphisms  $i_n$ and $p_n$ look like:
\[ i_n(\sum_{j=1}^l[c_j]\otimes_{M_0} \gamma_j)= \sum_{j=1}^l[c'_j]\otimes_{M_1}
\gamma_j+ \sum_{j=1}^l[c''_j]\otimes_{M_2} \gamma_j
\]
where $\gamma_j\in \mathbb{Z}M$, $j=1,2,\dots l$,
\[
\begin{array}{l}
p_n(\displaystyle\sum_{i=1}^m[a_i]\otimes_{M_1} \alpha_i +
\displaystyle\sum_{j=1}^l[c'_j]\otimes_{M_1} \delta'_j +
\displaystyle\sum_{i=1}^k[b_i]\otimes_{M_2} \beta_i +
\displaystyle\sum_{j=1}^l[c''_j]\otimes_{M_2}
\delta''_j)=\\
 =\displaystyle\sum_{i=1}^m[a_i]\alpha_i - \displaystyle\sum_{i=1}^k [b_i]\beta_i +
\displaystyle\sum_{j=1}^l[c_j](\delta'_j-\delta''_j).
\end{array}
\]

Evidently, for each $n$ the identity $p_ni_n=0$ holds. Besides,
the element $x= \sum_{i=1}^m[a_i]\otimes_{M_1} \alpha_i +
\sum_{j=1}^l[c'_j]\otimes_{M_1} \delta'_j +
\sum_{i=1}^k[b_i]\otimes_{M_2} \beta_i +
\sum_{j=1}^l[c''_j]\otimes_{M_2} \delta''_j$ belongs to ${\rm
Ker}\,p_n$ if and only if $\alpha_i=0$, $(i=1,\dots ,m)$,
$\beta_i=0$, $(i=1,\dots ,k)$ and $\delta'_j=\delta''_j$,
($j=1,\dots ,l$), i.e.
$x=i_n(\sum_{j=1}^l[c_j]\otimes_{M_0}\delta'_j)$. Hence, ${\rm
Im}\,i_n = {\rm Ker}\,p_n$ holds.

The exactness of all the rows, the left and the middle columns
implies the exactness of the right column what can be shown by
standard diagram search method. Thus the theorem is proved.
\rule{1ex}{1ex}

\begin{cor}\label{Cor1}
If  the graph $\Gamma(M)$ of free partially commutative monoid $M$
does not contain complete subgraphs with more than  $n$ vertices,
then the homological dimension of $M$ does not exceed $n$.
\end{cor}
Corollary~\ref{Cor1} proves the Husainov's Conjecture formulated
in the introduction.

\begin{cor}\label{Cor2}
Let $M$ be a free partially commutative monoid, $A$ be a trivial
left  $M$-module. Then for $n\geq 1$ homology groups
$H_n(M,A)\cong\underset{r_n}{\underbrace{A\oplus\dots \oplus A}}$
where $r_n$ is the number of complete subgraphs with $n$ vertices
in graph $\Gamma(M)$.
\end{cor}
{\bf Proof.} Denote $C_n = \underset{r_n}{\underbrace{A\oplus\dots
\oplus A}}$. Notice that there exists a homomorphism $F_n\otimes
A\cong C_n$. Indeed the mapping $\varphi: F_n\times A \rightarrow
C_n$ can be constructed in the following way. Let $[x_1],\dots
,[x_{r_n}]$ be the generators of $F_n$. For each $j=1,\dots ,r_n$,
$\alpha\in \mathbb{Z}M$, $a\in A$ we set $\varphi([x_j]\alpha,
a)=(0, \dots , \alpha a,\dots ,0)$, where $\alpha a$ is put on the
$j$-th place, and extend this mapping linearly. It is easy to
check that the Abelian group $C_n$ together with the mapping
$\varphi$ satisfies the universal property of the tensor product
$F_n\otimes_M A$ and, therefore, it is isomorphic to this tensor
product.

The groups $H_n(M,A)$ are the homology groups of the complex

\[\dots \rightarrow
C_3 \stackrel{\partial_3}\rightarrow
C_2\stackrel{\partial_2}\rightarrow C_1
\stackrel{\partial_1}\rightarrow  A
\]
where $\partial_n= \delta_n\otimes_M A$. Since $M$-module $A$ is
trivial, then
\[ \partial_k([a_1a_2\dots a_k]\otimes_M 1)=\sum_{j=1}^k[a_1\dots \widehat{a_j}\dots
a_k](a_j-1)(-1)^{j-1}\otimes_M 1=0.
\]
Hence,
\[H_n(M,A)={\rm Ker}\, \partial_n/{\rm Im}\, \partial_{n+1}\cong
C_n.
\]
\rule{1ex}{1ex}

The particular case of this Corollary is

\begin{cor}\label{Cor3}
Let $M$ be a free partially commutative monoid. Then the homology
groups $H_n(M,\mathbb{Z})$, $n\geq 1$ are the free Abelian groups
of rank $r_n$ where  $r_n$ is the number of complete subgraphs
with $n$ vertices in graph $\Gamma(M)$.
\end{cor}

\pagebreak

\end {document}